# Exploring 'pseudorandom' value addition operations in datasets:

## A layered approach to escape from normal-Gaussian patterns

*Ergon Cugler de Moraes Silva*

Getulio Vargas Foundation (FGV)
University of São Paulo (USP)
São Paulo, São Paulo, Brazil

contato@ergoncugler.com
www.ergoncugler.com

**Abstract**

In the realm of statistical exploration, the manipulation of pseudo-random values to discern their impact on data distribution presents a compelling avenue of inquiry. This article investigates the question: 'Is it possible to add pseudo-random values without compelling a shift towards a normal distribution?'. Employing Python techniques, the study explores the nuances of pseudo-random value addition within the context of additions, aiming to unravel the interplay between randomness and resulting statistical characteristics. The Materials and Methods chapter details the construction of datasets comprising up to 300 billion pseudo-random values, employing three distinct layers of manipulation. The Results chapter visually and quantitatively explores the generated datasets, emphasizing distribution and standard deviation metrics. The study concludes with reflections on the implications of pseudo-random value manipulation and suggests avenues for future research. In the layered exploration, the first layer introduces subtle normalization with increasing summations, while the second layer enhances normality. The third layer disrupts typical distribution patterns, leaning towards randomness despite pseudo-random value summation. Standard deviation patterns across layers further illuminate the dynamic interplay of pseudo-random operations on statistical characteristics. While not aiming to disrupt academic norms, this work modestly contributes insights into data distribution complexities. Future studies are encouraged to delve deeper into the implications of data manipulation on statistical outcomes, extending the understanding of pseudo-random operations in diverse contexts.

## 1. Introduction

In the landscape of statistical exploration, the intricacies of manipulating pseudo-random values to discern their impact on data distribution present a compelling avenue of inquiry. In his seminal work published in the 90s, James (1990) embarked on a comprehensive review of pseudorandom number generators, shedding light on a critical aspect of computational processes. His examination revealed that while pseudorandom number generators with essential properties were accessible, the generators commonly in use often fell short of meeting the required standards.

Building on this foundation, Lagarias (1993) expanded the discourse by introducing the theory of computational information. At the heart of this theory lies the quest for a pseudorandom bit generator that not only ensures randomness but also guarantees security without introducing biases. Intriguingly, despite the apparent necessity and practical demand for such generators, their existence had not been conclusively proven through theoretical means at the time. Nonetheless, various functions were identified that exhibited the desired properties, sparking optimism within the computational community.

It is crucial to note that Lagarias' exploration (1993) did not merely stop at theoretical propositions; rather, it extended into the realm of practical applicability. Despite the lack of a robust theoretical proof for the existence of perfectly secure pseudorandom generators, empirical evidence suggested that certain generators, while not theoretically proven, demonstrated commendable performance in real-world scenarios.

In the broader landscape of computational science, the interplay between theory and practice in the realm of pseudorandom number generators poses intriguing questions. While theoretical developments provide a roadmap and a conceptual framework for understanding the properties of ideal generators, the practical efficacy of certain generators cannot be denied. The field stands at the intersection of theoretical aspirations and pragmatic demands, where ongoing research seeks to bridge the gap between theory and practice in the pursuit of reliable and secure pseudorandom number generation.

In recent publications, the implementation of pseudorandomness extends across various domains, from functions within the Internet of Things (IoT) (Kietzmann et al., 2021) to the utilization of pseudo-random states as a foundation for cryptographic tasks (Ananth et al., 2022). Addressing the secure generation of correlated randomness with minimized communication costs, Boyle et al. (2022) propose a groundbreaking design for a pseudorandom correlation generator (PCG).

This design, grounded in expand-accumulate codes, not only ensures competitive concrete efficiency with provable security but also introduces an offline-online mode to optimize cache-friendliness and parallelization. It provides extensions to pseudorandom correlation functions, enabling the incremental generation of new correlation instances. The authors also present a method to enhance the evaluation of puncturable pseudorandom functions (PPRF), a concept independently motivated by various cryptographic applications.

Concurrently, Dwork et al. (2023) delve into the intricate connections between multi-group fairness in prediction algorithms and pseudorandomness concepts. Their exploration, employing statistical distance-based variants of multicalibration, not only yields graph theoretic results but also leads to the development of more efficient multicalibration algorithms. Additionally, their work presents a novel proof of a hardcore lemma for real-valued functions, contributing to the theoretical underpinnings of pseudorandomness in prediction algorithms. This interconnected research landscape reflects the growing importance of pseudorandomness across diverse applications, showcasing its relevance in the realms of IoT, cryptography, and prediction algorithms.

This integrated narrative sets the stage for the current article's exploratory journey, probing the nuanced question of adding pseudo-random values without inducing a shift towards a normal distribution. In this sense, this article embarks on an exploratory journey, delving into the nuanced question: **'Is it possible to add pseudo-random values without compelling a shift towards a normal distribution?'**. This query encapsulates the essence of an investigation that seeks not only to unravel the complexities of pseudo-random value addition but also to shed light on the subtle interplay between randomness and the resulting statistical characteristics.

Motivated by the need to comprehend the underlying patterns governing data distribution, this study employs Python techniques as its investigative toolbox. The chosen approach involves encapsulating pseudo-randomness within the realm of additions, striving to discern how this manipulation influences the resultant distribution patterns. The objective is not to disrupt established academic norms but rather to provide a nuanced understanding of a phenomenon that has implications across various disciplines.

The narrative unfolds through distinct chapters, each contributing a unique facet to the overarching exploration. The Materials and Methods chapter intricately details the construction of datasets comprising up to 300 billion pseudo-random values. This substantial dataset serves as the canvas for the application of three distinct layers, each representing a novel approach to manipulating pseudo-random values. The Python implementation of these layers stands as a testament to the rigor and precision applied in this exploration.

Moving forward, the Results chapter immerses the reader in a visual and quantitative exploration of the generated datasets. Distribution and standard deviation metrics take center stage, offering a detailed panorama of the evolving statistical landscape as pseudo-random values undergo varying degrees of manipulation. Through detailed visualizations and metrics, this section unravels the subtle shifts and transformations within the datasets, highlighting the impact of each layer on the distribution patterns.

As the exploration culminates, the concluding chapter unfolds a tapestry of reflections and challenges, paving the way for future investigations. Far from presenting conclusive answers, this work prompts contemplation on the intricacies of pseudo-random value manipulation, urging researchers to consider auxiliary techniques in pseudo-randomization within the context of summing pseudo-random values.

## 2. Materials and methods

Using Python techniques, 10,000 datasets are generated, each containing 1,000 numbers. Each number is the result of an addition of up to 10,000 random numbers. Considering the upper limit of randomly generated values to yield a addition product, we have a scenario of up to 100,000,000,000 (one hundred billion) generated values. As tests are conducted across three layers, these values are tripled, reaching an estimated total of 300,000,000,000 (three hundred billion) generated values.

## 2.1. First layer

The approach involves the absence of added randomness within the set or quantity loop. The process unfolds with progressive additions as the loop iterates. Within this layer, the data is systematically generated, commencing with the first set containing a single random number. Subsequent sets build upon this progression, where each set incorporates an additional random number. This cumulative pattern continues until the last set, which culminates in the addition of 10,000 random numbers.

The nested loops guide this structured generation, with the outermost loop iterating over sets, the middle loop handling the quantity of numbers in each set, and the innermost loop orchestrating the random additions. In this sense, Figure 01, below, illustrates how the distribution of random values generated for each dataset occurs within a layer.

**Figure 01. Pseudo-random values**

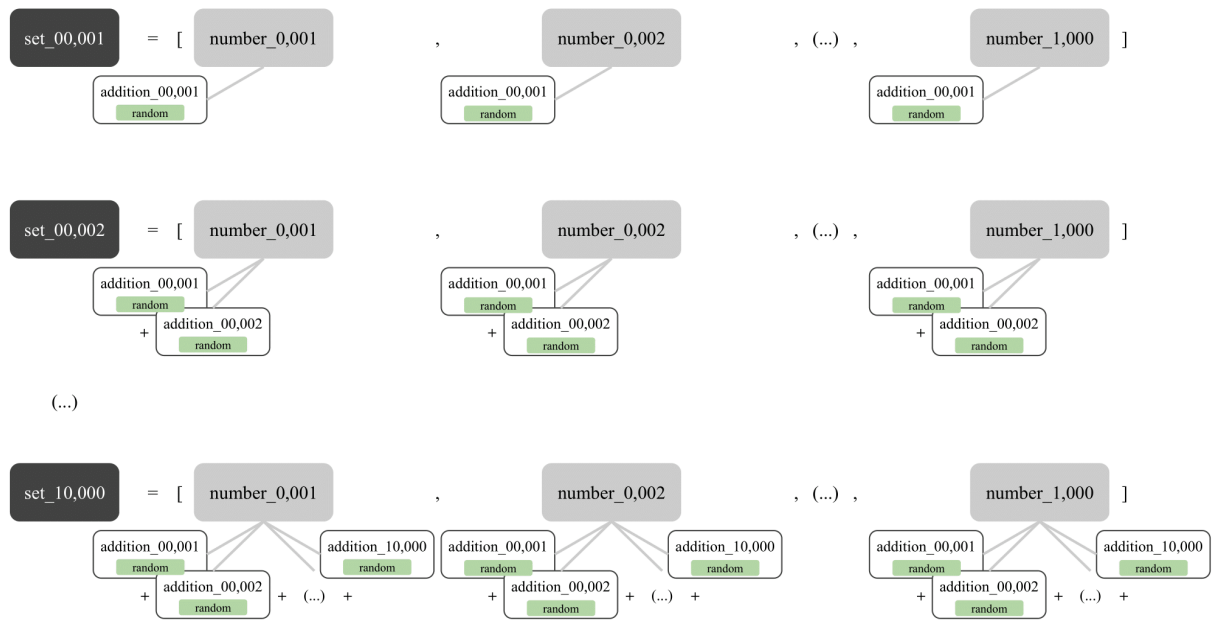

**Source:** Own elaboration (2024).

$$v_{i,j,k} = v_{i-1,j,k} + random\,(1,\ max\_number)$$

Where $v_{i,j,k}$ is a random number between 1 and max_number, $i$ ranges from 2 to total_additions, $j$ ranges from 1 to total_numbers, e $k$ ranges from 1 to total_sets.

## 2.2. Second layer

In the 2nd Pseudo-Randomity Layer, the technique introduces randomness inside the set loop while maintaining a consistent amount across the quantity loop. Each number within a set undergoes the same random amount of additions, ranging from 1 to 10,000. This variation leads to diverse sets, each with a distinct total addition. The implementation involves nested loops, where the outer loop controls the sets, the middle loop manages the quantity of numbers, and the innermost loop orchestrates the random additions. This layer allows for variability in the total addition of each set, offering a nuanced perspective on the impact of random additions on the generated data (Figure 02).

**Figure 02. Sets with pseudo-random additions**

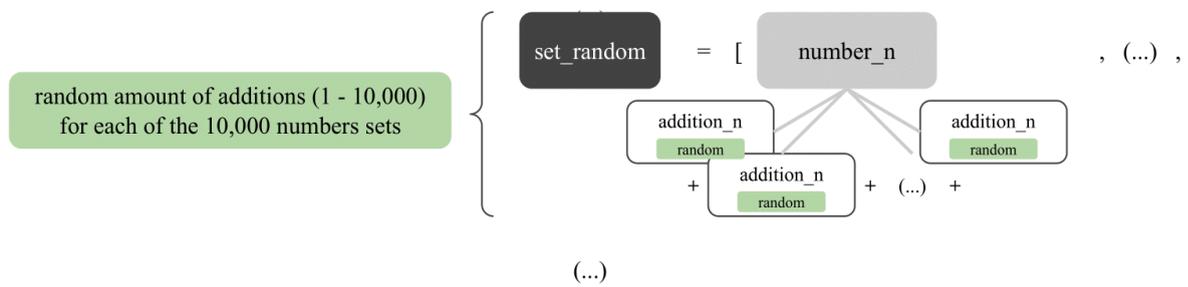

**Source:** Own elaboration (2024).

$$random\_additions_k = random(1, total\_additions)$$

$$v_{i,j,k} = v_{i-1,j,k} + random(1, max\_number)$$

Where $v_{i,j,k}$ is a random number between 1 and max_number, $i$ ranges from 2 to random_additions$_k$, $j$ ranges from 1 to total_numbers, e $k$ ranges from 1 to total_sets.

## 2.3. Third layer

The 3rd Pseudo-Randomity Layer adopts an approach where randomness is incorporated within the quantity loop, creating a dynamic generation process for each set (Figure 03). In this layer, each number within a set possesses a random quantity of additions. The intricacy of this layer lies in the flexibility it provides, allowing each of the 1,000 numbers in a set to have a distinct number of random additions. The randomness ranges from 1 to 10,000, offering a wide spectrum of total additions within a set. The nested loops govern this intricate process, with the outer loop iterating over sets, the middle loop managing the quantity of numbers, and the innermost loop orchestrating the random additions. This layer introduces a granular level of variability in the additions generated, contributing to a comprehensive understanding of the impact of random additions on the resulting dataset.

**Figure 03. Each number within the sets with number of pseudo-random additions**

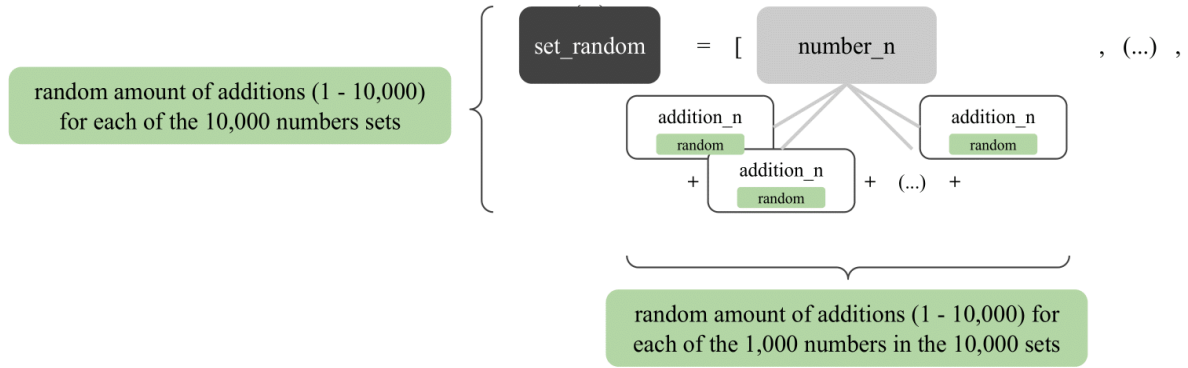

**Source:** Own elaboration (2024).

$$random\_additions_{j,k} = random(1, total\_additions)$$
$$v_{i,j,k} = v_{i-1,j,k} + random(1, max\_number)$$

Where $v_{i,j,k}$ is a random number between 1 and max_number, $i$ ranges from 2 to random_additions$_{j,k}$, $j$ ranges from 1 to total_numbers, e $k$ ranges from 1 to total_sets.

### 2.4. Codes

We employed a variety of Python libraries to conduct our statistical analyses and visualization tasks. The **scipy.stats** library provided essential statistical functions such as skewness, kurtosis, Shapiro-Wilk test, Kolmogorov-Smirnov test, and the normal distribution. For streamlined and interactive data visualization, we leveraged the capabilities of **plotly.express**, enabling us to create insightful plots and charts. Data manipulation and handling were efficiently managed using the versatile **pandas** library, while numerical operations and array handling were facilitated by **numpy**. Additionally, the **statistics** library proved useful for basic statistical calculations, and the **random** module allowed for the generation of random numbers when necessary. Lastly, we utilized the **time** module to measure and monitor certain processes. In Table 01, the variables are set as **max_number** = 100; **total_numbers** = 1000; **total_sets** = 10000; and **total_additions** = 10000.

Table 01. Approaches and codes

| Approach description | Code description |
|---|---|
| **1st Pseudo-Randomity Layer:**<br><br>No added randomness within the set or quantity loop, just progressive additions as you go through the loop. Here the data is generated with additions progressively. In the first set, there is just a random number; In the second set, there are two random numbers added together; (...) In the last set, there are 10,000 random numbers added together. | ```python<br># Sets loop<br>for set in range(1, total_sets + 1):<br>    data = []<br><br>    # Amount loop<br>    for _ in range(total_numbers):<br>        value = random.randint(1, max_number)<br><br>        # Addition loop<br>        for loop in range(total_additions - 1):<br>            value = value + random.randint(1, max_number)<br>        data.append(value)<br>``` |
| **2nd Pseudo-Randomity Layer:**<br><br>Addition of randomness inside the set loop but outside the quantity loop, i.e. each number in the set contains the same random amount of additions. Here data is generated with random additions. The first set can have between 1 and 10,000 random numbers added together; The second set can have between 1 and 10,000 random numbers added together; (...) The last set can have between 1 and 10,000 random numbers added together. | ```python<br># Sets loop<br>for set in range(1, total_sets + 1):<br>    data = []<br>    random_additions = random.randint(1, total_additions)<br><br>    # Amount loop<br>    for _ in range(total_numbers):<br>        value = random.randint(1, max_number)<br><br>        # Addition loop<br>        for loop in range(random_additions - 1):<br>            value = value + random.randint(1, max_number)<br>        data.append(value)<br>``` |
| **3rd Pseudo-Randomity Layer:**<br><br>Addition of randomness within the quantity loop, i.e. each number in the set contains a random quantity of additions. Here the data is generated with random additions within the sets themselves. Each of the 1,000 numbers in the first set can have between 1 and 10,000 random numbers added together; Each of the 1,000 numbers in the second set can have between 1 and 10,000 random numbers added together; (...) Each of the 1,000 numbers in the last set can have between 1 and 10,000 random numbers added together. | ```python<br># Sets loop<br>for set in range(1, total_sets + 1):<br>    data = []<br><br>    # Amount loop<br>    for _ in range(total_numbers):<br>        random_additions = random.randint(1, total_additions)<br>        value = random.randint(1, max_number)<br><br>        # Addition loop<br>        for loop in range(random_additions - 1):<br>            value = value + random.randint(1, max_number)<br>        data.append(value)<br>``` |

**Source:** Own elaboration (2024).

## 3. Exploring the results

The following Figure 04 visually depicts the distribution of pseudo-random values before any summation. As indicated by the literature, this distribution exhibits characteristics of a 'non-normal' or 'non-Gaussian' distribution. The variability in the spread of values along the scale is evident, emphasizing the absence of a typical symmetrical shape associated with a normal distribution. This initial pattern provides a foundation for understanding the effect of summation on the distribution's properties (Figure 04).

**Figure 04. Distribution of random values without addition**

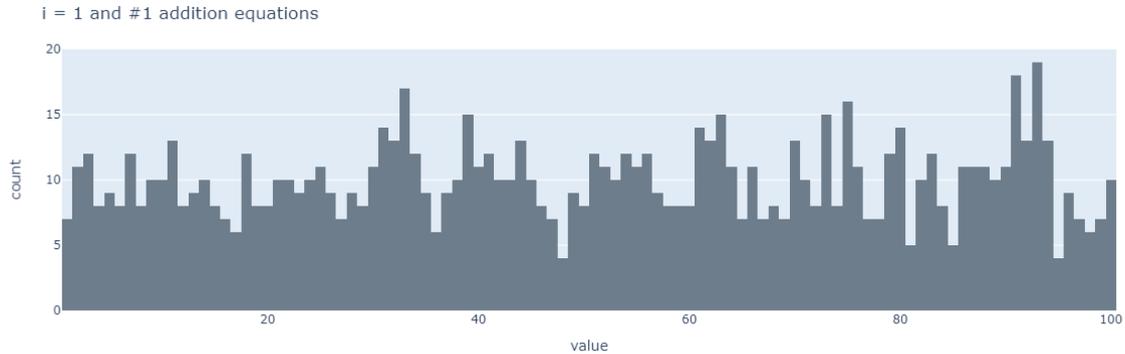

**Source:** Own elaboration (2024).

However, upon adding a single summation to each value in the dataset (pseudo-random value + pseudo-random value), a striking transformation in the distribution's shape is observed, as seen in Figure 05. The representation reveals a triangular shape emerging, indicating the initiation of the formation of a normal or Gaussian distribution. This phenomenon highlights the significant influence of the summation operation on transforming the original distribution, showcasing the impact of data manipulation on its statistical shape.

**Figure 05. Distribution of random values with an addition**

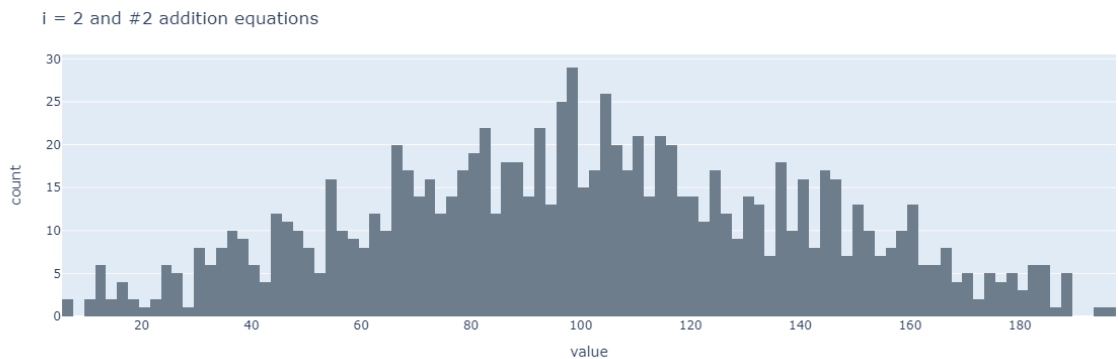

**Source:** Own elaboration (2024).

In the expanded example in Figure 06, we delve deeper into the curve normalization phenomenon by performing 10 summations of pseudo-random variables. The distribution, which initially displayed a triangular shape, now exhibits a clear Gaussian trend. The progressive smoothing of the curve and the broadening of the base indicate convergence toward a normal distribution as the number of summations increases. This visualization reinforces the observation that the successive summation of pseudo-random values tends to shape the distribution towards a more Gaussian form.

**Figure 06. Distribution of random values with multiple additions**

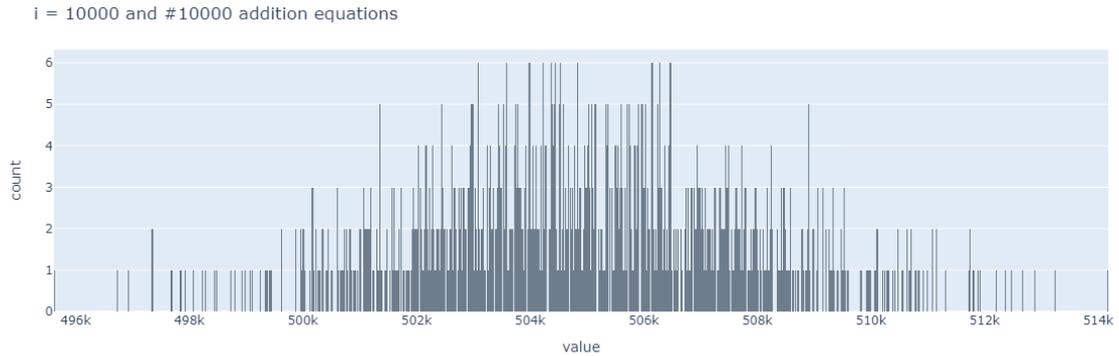

**Source:** Own elaboration (2024).

Having delved into the fundamental aspects of pseudo-random value addition operations, the subsequent sections provide a comprehensive exploration of the layer implementations designed to counteract the tendency of normalizing data distribution.

### 3.1. Distribution with layers

In the first layer, a nuanced examination of the data distributions is conducted by scaling boxplots to an equivalent scale. Figure 07 intricately portrays twelve distribution examples, ranging from the initial boxplot with no additions (addition = 1) to those with up to 10,000 additions. This comprehensive visualization captures the evolving nature of the distribution pattern as the number of additions increases.

**Figure 07. Boxplots in first layer implementation**

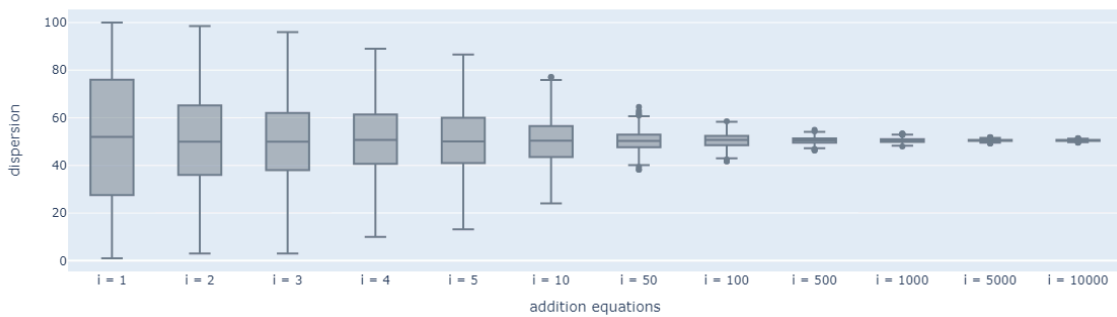

**Source:** Own elaboration (2024).

Upon scaling in terms of logical sets encompassing 10,000 datasets, a shift in perspective is observed. Rather than scrutinizing individual datasets in isolation, the addition of a second layer reveals intriguing insights. This layer not only curtails outliers by defining the quantity of additions to be performed but also augments the normalization of the data, as depicted in Figure 08. The interplay between layers becomes more apparent when viewed from this broader perspective, unveiling the intricate dynamics shaping the data distribution.

**Figure 08. Boxplots in second layer implementation**

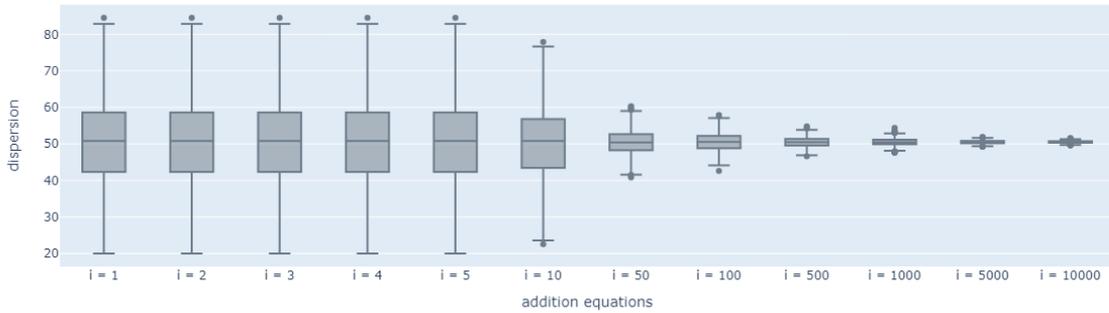

**Source:** Own elaboration (2024).

However, the introduction of the third layer of pseudo-randomness, determining the quantity of additions each value will undergo within the datasets, disrupts the normal distribution characteristics of the data. Figure 09 provides a visual representation of this aspect, showcasing boxplots with an almost identical distribution and signaling a departure from the Gaussian pattern. This intricate exploration sheds light on the intricate interplay of pseudo-random operations and their impact on the statistical properties of the data.

**Figure 09. Boxplots in third layer implementation**

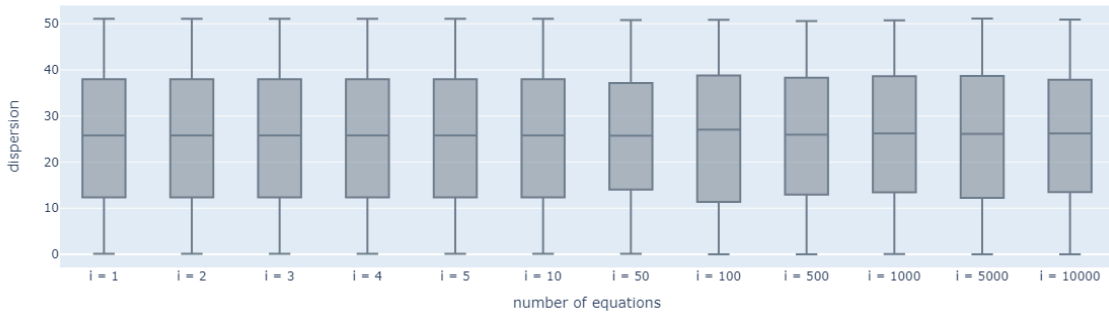

**Source:** Own elaboration (2024).

In alignment with the insights provided by the boxplot analysis, a more nuanced exploration of the data distribution is afforded by Figures 10 and 11 below. Figure 10 encapsulates the impact of the first layer, showcasing how the initial implementation of pseudo-random value addition introduces a degree of normalization. The boxplots artfully depict the evolving distribution pattern, emphasizing the initial stages of Gaussian tendencies as randomness is incorporated into the data. Building upon this, Figure 11 delves into the second layer, where the limitation of outliers and the controlled definition of the quantity of additions further contribute to the development of a more normalized data distribution.

**Figure 10. Distribution in first layer implementation**

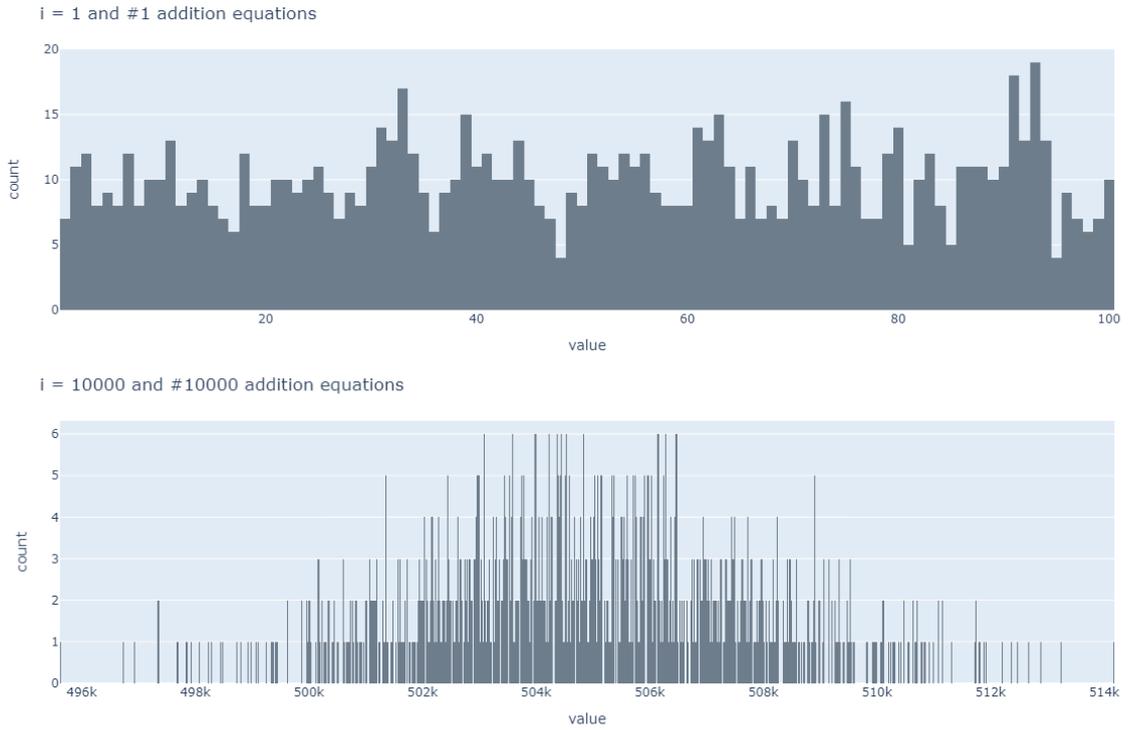

**Source:** Own elaboration (2024).

**Figure 11. Distribution in second layer implementation**

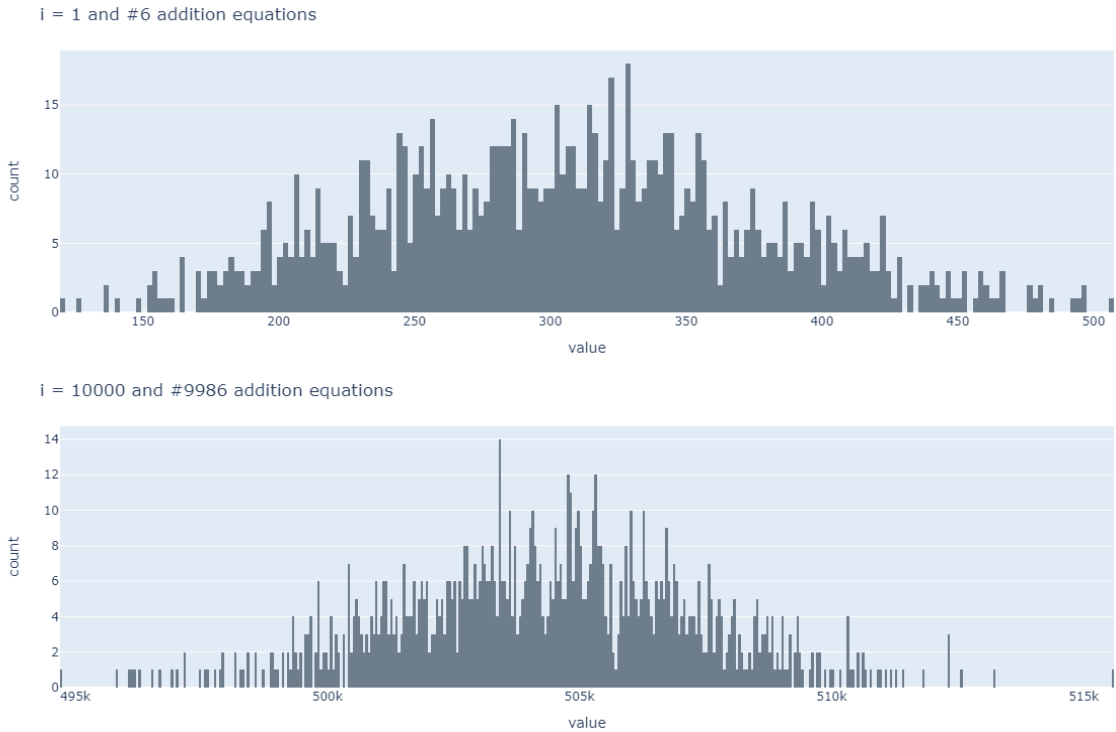

**Source:** Own elaboration (2024).

Yet, it is Figure 12 that sheds light on the phenomenon at hand, where the implementation of the third layer in the summation of pseudo-random values made it possible to distribute data without adhering to a Gaussian pattern. This figure serves as a pivotal insight into how the nuanced manipulation of random quantities in the addition operations disrupts the typical trend towards normalization. In Figure 12, the intricacies of the third layer's impact become evident as it introduces randomness into the quantity of addition operations. The boxplots illustrate the divergence from a Gaussian pattern, showcasing the distribution's resilience against normalization biases. The inclusion of this additional layer further underscores the effectiveness of randomization in countering tendencies towards standardization, presenting a nuanced perspective on how pseudo-random operations shape the statistical characteristics of the data.

**Figure 12. Distribution in third layer implementation**

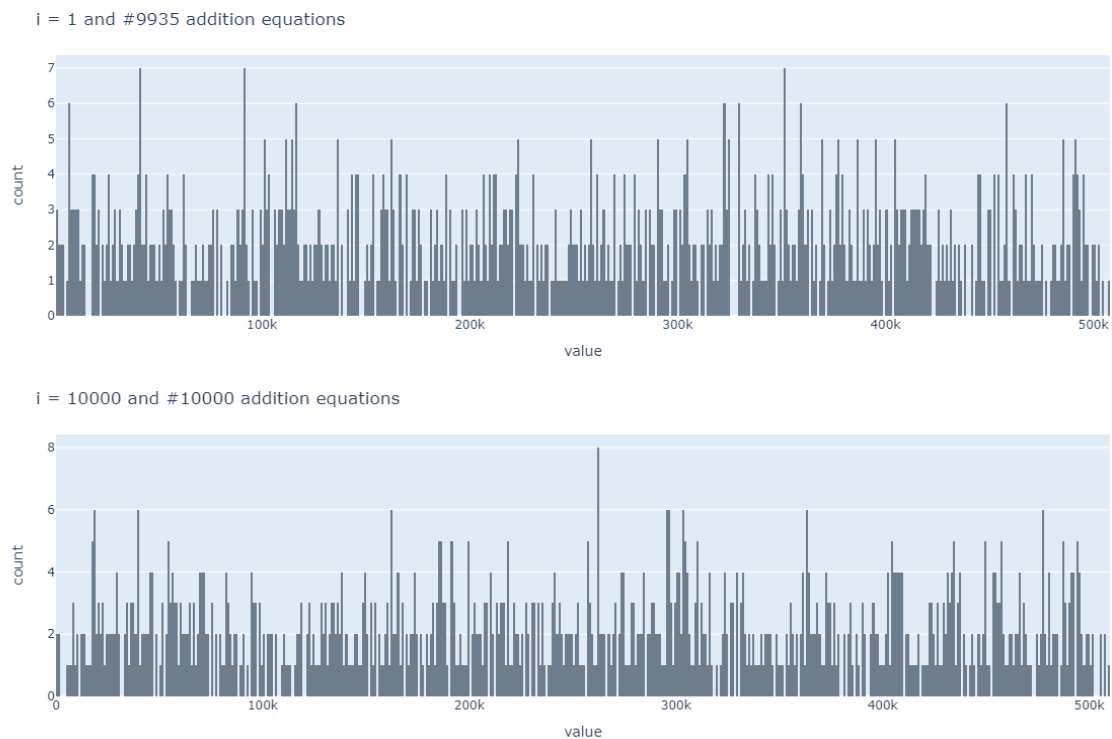

**Source:** Own elaboration (2024).

This investigation into distribution patterns employing layered pseudo-random operations unravels a dynamic interplay between randomness and the statistical characteristics of the data. Beginning with the initial layer, Figure 07 portrays the evolving distribution pattern as additions increase, introducing a subtle degree of normalization. Scaling up to logical sets of 10,000 datasets broadens our perspective, highlighting the second layer's role in curbing outliers and contributing to a more normalized distribution, as seen in Figure 08. However, the introduction of the third layer disrupts the conventional Gaussian trend, as illustrated in Figure 09, showcasing the potential to distribute data without adhering to typical patterns. Figures 10 and 11 offer a deeper understanding of the layered approach, depicting progressive normalization resulting from the first and second layers. Yet, Figure 12 stands out, showcasing how the third layer enables data distribution without conforming to a Gaussian pattern. This final figure emphasizes the pivotal role of nuanced randomization in the quantity of addition operations, highlighting its effectiveness in countering normalization biases.

### 3.2. Standard deviation with layers

The same dynamic appears to unfold when we examine the standard deviation. The datasets tested with the first layer exhibit a distribution that ranges from approximately 30% (of the total values) to less than 1% (Figure 13). This variability suggests a considerable dispersion within the dataset. About the second layer (Figure 14), we observe a more concentrated distribution as the randomization of sets minimizes the presence of extreme values. However, t is the datasets tested with a third layer (Figure 15) that reveal a non-Gaussian distribution, characterized by an almost linear pattern of standard deviation percentages. Upon closer inspection of Figure 13, the diversity in the standard deviation percentages for datasets subjected to the first layer becomes evident. This wide-ranging distribution implies substantial variability in the spread of values within the dataset.

**Figure 13. Distribution of standard deviation in first layer implementation**

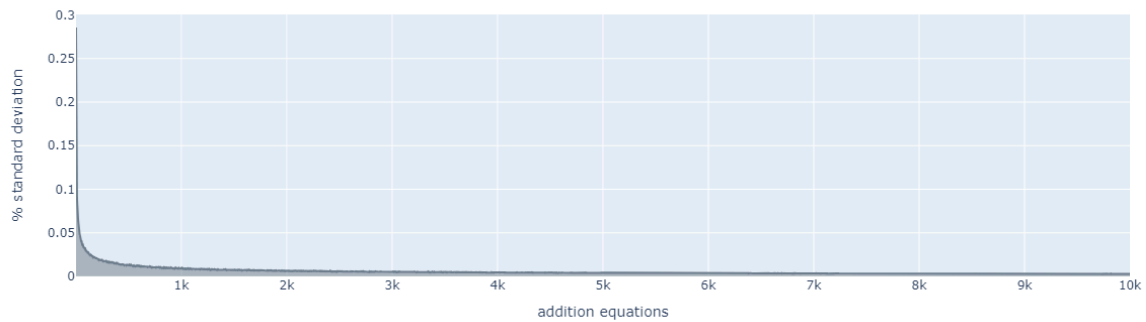

**Source:** Own elaboration (2024).

In Figure 14, focusing on datasets subjected to the second layer, we witness a more focused distribution. The introduction of set randomization contributes to a reduction in extreme values, resulting in a narrower range of standard deviation percentages.

**Figure 14. Distribution of standard deviation in second layer implementation**

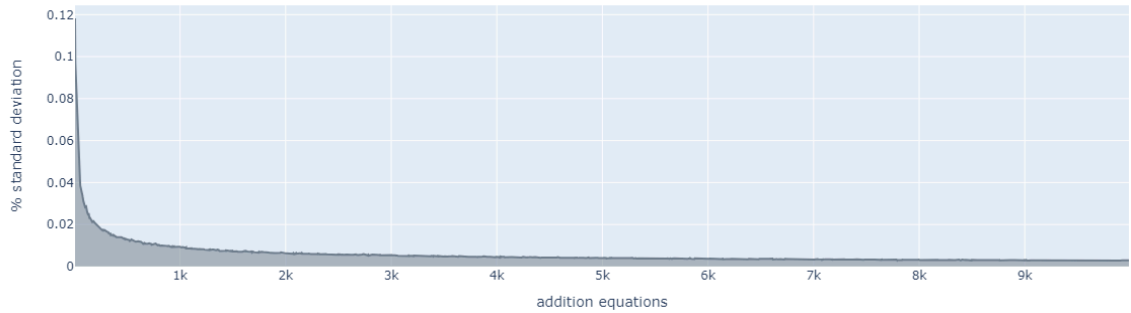

**Source:** Own elaboration (2024).

Figure 15 presents a compelling snapshot of datasets subjected to the third layer, showcasing a distribution with a distinct departure from Gaussian patterns. The almost linear pattern in the percentage of standard deviation emphasizes the unique impact of the third layer. This layer, which introduces randomness into the quantity of addition operations, appears to drive the standard deviation percentages towards a more evenly spaced and non-Gaussian distribution. This detailed exploration of standard deviation patterns across layers provides additional insights into the intricate dynamics of pseudo-random operations and their influence on statistical characteristics.

**Figure 15. Distribution of standard deviation in third layer implementation**

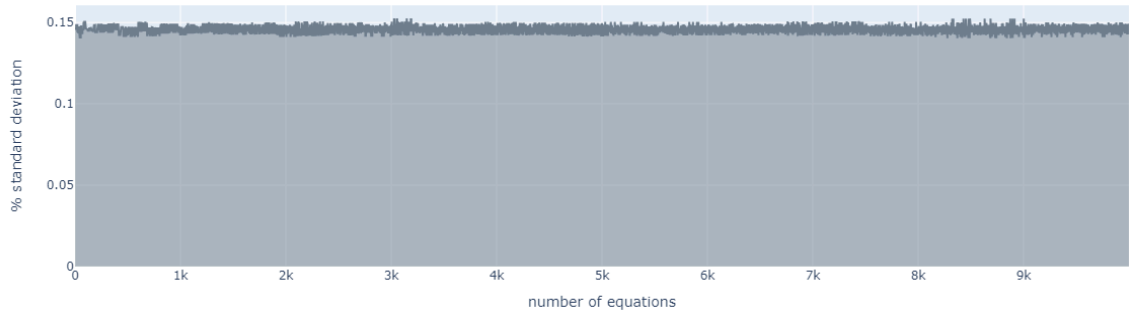

**Source:** Own elaboration (2024).

In conclusion, the examination of standard deviation patterns across layers provides intriguing insights into the dynamic interplay of pseudo-random operations and their influence on statistical characteristics. The datasets subjected to the first layer exhibit a wide-ranging distribution, indicating substantial variability within the dataset. Transitioning to the second layer results in a more focused distribution as set randomization minimizes extreme values, narrowing the range of standard deviation percentages (Figure 14). However, it is the datasets tested with a third layer that reveal a non-Gaussian distribution, characterized by an almost linear pattern of standard deviation percentages (Figure 15).

4. **Statistical summarization**

As we can see in Table 02, below, although the mean and median remain in the center of the distribution in the first and second layers, it is in the third layer that there is a substantial change in the displacement of the distribution axis. Also, the standard deviation stands out, deviating from a normal distribution when we see the third layer.

**Table 02. Metrics in first layer implementation**

|  | Add | Average | % Average | Median | % Median | Std dev | % Std dev |
|---|---|---|---|---|---|---|---|
| **1st layer** | 00,001' round | 51 | 51,2% | 52 | 52,0% | 28.555.402 | 28,6% |
|  | (...) | | | | | | |
|  | 10,000' round | 504.982 | 50,5% | 504.863 | 50,5% | 2.827.526.251 | 0,3% |
| **2nd layer** | 00,001' round | 305 | 50,9% | 305 | 50,8% | 70.942.868 | 11,8% |
|  | (...) | | | | | | |
|  | 10,000' round | 504.339 | 50,5% | 504.431 | 50,5% | 2.948.602.495 | 0,3% |
| **3rd layer** | 00,001' round | 252.118 | 25,4% | 256.306 | 25,8% | 145.169.707.627 | 14,6% |
|  | (...) | | | | | | |
|  | 10,000' round | 256.254 | 25,6% | 261.975 | 26,2% | 145.334.810.222 | 14,5% |

**Source:** Own elaboration (2024).

5. **Reflections and future work**

This article aimed to investigate whether **'Is it possible to add pseudo-random values without tending to a normal distribution?'**. To address this issue, techniques in Python were employed for hypothesis testing with pseudo-randomization in 'three layers'. Through the generation of 10,000 datasets, each containing 1,000 numbers, and employing up to 10,000 random additions for each number, a vast and diverse dataset of 300,000,000,000 (three hundred billion) values was created. The layers were meticulously designed to explore the impact of random additions on data distribution.

The first layer, marked by progressive additions within the set and quantity loops, introduced a subtle degree of normalization. Scaling up to the second layer, randomness was injected within the set loop, curbing outliers and contributing to a more normalized distribution. However, it was the third layer, incorporating randomness within the quantity loop, that disrupted the conventional Gaussian trend, offering the potential to distribute data without adhering to typical patterns. The layered approach illustrated in Figures 10 to 11 highlights the nuanced role of randomization in countering normalization biases.

Exploring the results further, the analysis extended to standard deviation patterns across layers. The datasets subjected to the first layer exhibited a wide-ranging distribution, indicating substantial variability. The second layer, focusing on set randomization, resulted in a more concentrated distribution by minimizing extreme values. Nevertheless, it was the datasets tested with the third layer that revealed a non-Gaussian distribution, characterized by an almost linear pattern of standard deviation percentages. Therefore, Table 03 summarizes:

**Table 03. Summary of findings**

| Layer | What is it | What was found |
|---|---|---|
| 1st | Involves the summation of pseudo-random numbers. As the quantity of additions increases, the distribution tends to become more normal (a concept well-established in existing literature). | The distribution becomes progressively more normal as the number of summations increases, aligning with existing literature. |
| 2nd | Encompasses the summation of pseudo-random numbers, with each of the 10,000 sets varying in the range of 1 to 10,000 additions. The distribution of values with this second layer results in an even more normal (bell-shaped) curve compared to a single layer. | The second layer contributes to a distribution with enhanced normality, surpassing the effects of a single layer. |
| 3rd | Involves the summation of pseudo-random numbers, with each of the 1,000 numbers within the 10,000 sets varying in the range of 1 to 10,000 additions. The distribution of values with this third layer results in a lack of a non-normal distribution, tending towards randomness, even with the summation of pseudo-random values. | The third layer leads to a distribution devoid of non-normal patterns, leaning towards randomness despite the addition of pseudo-random values. |

**Source:** Own elaboration (2024).

In conclusion, the exploration of distribution patterns and standard deviation across layers illuminates the dynamic interplay of pseudo-random operations on statistical characteristics. While the first and second layers introduced progressive normalization, it was the third layer that showcased the potential to disrupt Gaussian trends. This layered approach represents a modest exploration of data distribution, aiming to understand the potential impact of pseudo-random value manipulation on statistical properties.

While our findings contribute to the ongoing conversation surrounding pseudo-random values and their influence on data distribution, it is essential to acknowledge the limitations of this study. Our work does not aspire to introduce groundbreaking academic disruptions but rather to utilize available techniques for investigation. Future studies are warranted to delve deeper into the complexities of data manipulation, offering a more comprehensive understanding of its implications on statistical outcomes.

## 7. Author biography

**Ergon Cugler de Moraes Silva** has a Master's degree in Public Administration and Government (FGV), Postgraduate MBA in Data Science & Analytics (USP) and Bachelor's degree in Public Policy Management (USP). He is associated with the Bureaucracy Studies Center (NEB), collaborates with the Interdisciplinary Observatory of Public Policies (OIPP), with the Study Group on Technology and Innovations in Public Management (GETIP) and with the Monitor of Political Debate in the Digital Environment (Monitor USP). São Paulo, São Paulo, Brazil. Web site: https://ergoncugler.com/.